\documentclass[a4paper,twoside,10pt,reqno]{amsart}
\usepackage{amsmath,amsfonts,amssymb,amsthm}
\usepackage{mathrsfs}
\usepackage{graphicx}
\usepackage{color}
\usepackage[a4paper]{geometry}
\usepackage{cite}
\usepackage[colorlinks=true,pdfstartview={XYZ null null 1.00},pdfview=FitH,citecolor=blue,urlcolor=customgreen]{hyperref}
\newtheorem{theorem}{Theorem}
\numberwithin{equation}{section}
\def\im{\mathrm{i}}
\def\d{\mathrm{d}}
\def\id{\,\mathrm{d}}
\def\e{\mathrm{e}}

\def\Re{\operatorname{Re}}

\def\Im{\operatorname{Im}}
\definecolor{customgreen}{rgb}{0.0, 0.5, 0.0}

\def\Yint#1{\mathchoice
    {\YYint\displaystyle\textstyle{#1}}%
    {\YYint\textstyle\scriptstyle{#1}}%
    {\YYint\scriptstyle\scriptscriptstyle{#1}}%
    {\YYint\scriptscriptstyle\scriptscriptstyle{#1}}%
      \!\int}
\def\YYint#1#2#3{{\setbox0=\hbox{$#1{#2#3}{\int}$}
    \vcenter{\hbox{$#2#3$}}\kern-.52\wd0}}
\def\cupint{\:\,\Yint\cup}

\author[D. Broadhurst]{David Broadhurst}
\address{School of Physical Sciences, The Open University, Walton Hall, Kents Hill, Milton Keynes MK7 6AA, United Kingdom}

\author[G. Nemes]{Gerg\H{o} Nemes}
\address{School of Mathematics, Harbin Institute of Technology, Xidazhi Street, Harbin 150001,\linebreak Heilongjiang, People's Republic of China\medskip}

\email{djbroadhurst@yahoo.co.uk}
\email{g.nemes@hit.edu.cn}


\begin{document}

\title[Polynomials and asymptotic constants]{Polynomials and asymptotic constants\\
in a resurgent problem from {}'t Hooft}

\begin{abstract}
In a recent study of the quantum theory of harmonic oscillators, Gerard {}'t Hooft proposed the following problem: given $G(z)=\sum_{n=1}^\infty\sqrt{n}\,z^n$ for $|z|<1$, find its analytic continuation for $|z|\ge1$, excluding a branch-cut $z\in[1,\,\infty)$. A solution is provided by the bilateral convergent sum
$G(z)=\frac12\sqrt{\pi}\sum_{n=-\infty}^\infty(2\pi\im n-\log(z))^{-3/2}$. On the negative real axis, $G(-\e^u)$ has a sign-constant asymptotic expansion in $1/u^2$, for large positive $u$. Optimal truncation leaves exponentially suppressed terms in an asymptotic expansion $\e^{-u}\sum_{k=0}^\infty P_k(x)/u^k$, with $P_0(x)=x-\frac23$ and $P_k(x)$ of degree $2k+1$ evaluated at $x=u/2-\lfloor u/2\rfloor$. At large $k$, these polynomials become excellent approximations to sinusoids. The amplitude of $P_k(x)$ increases factorially with $k$ and its phase increases linearly, with $P_k(x)\sim\sin((2k+1)C-2\pi x) R^{2k+1}\Gamma(k+\frac12)/\sqrt{2\pi}$,
where $C\approx1.0688539158679530121571$ and $R\approx0.5181839789815558726739$ are asymptotic constants satisfying $R\exp(\im\,C)=\sqrt{-1/(2+\pi\im)}$.
\end{abstract}

\maketitle

\section{Introduction}

In~\cite{GH}, Gerard {}'t Hooft sought an analytic continuation of the fractional polylogarithm $\operatorname{Li}_{-1/2}(z)$, denoted there and herein by $G(z)$,
\begin{equation}
G(z)=\sum_{n=1}^\infty\sqrt{n}\,z^n,\text{ for } |z|<1\label{Gz}
\end{equation}
to the region $|z|\ge 1$, excluding the branch cut along $[1,\infty)$.

A solution to this problem is known~\cite{OC,DW}, namely
\begin{equation}
G(z)=\frac{\sqrt{\pi}}{2}\sum_{n=-\infty}^{\infty}\left(\frac{1}{2\pi\im n-\log(z)}\right)^{3/2},
\text{ for }z\notin[1,\infty).\label{Gsum}
\end{equation}

The paper is organised as follows. Section~\ref{section2} reviews computational strategies for $G(z)$. Section~\ref{section3} considers 
its behaviour on the negative axis where $G(-\e^u)$, at large positive $u$, is given
by optimal truncation of a sign-constant asymptotic expansion in $1/u^2$, leaving
an exponentially suppressed term $\e^{-u}S(u)$. Then the asymptotic expansion  of 
$S(u)\sim\sum_{k=0}^\infty P_k(x)/u^k$ yields a remarkable sequence of polynomials, $P_k(x)$,
with argument $x=u/2-\lfloor u/2\rfloor$.   Section~\ref{section4} gives formula \eqref{Pkasy},
for $P_k(x)$ at large $k$, and Theorem~\ref{thm1} identifies the constants $C$ and $R$ in this formula as solutions to 
$R\exp(\im\,C)=\sqrt{-1/(2+\pi\im)}$, as proved in Section~\ref{section5}. The paper concludes with a brief discussion in Section~\ref{section6}.

We note that while the general resurgence and monodromy of the polylogarithms $\operatorname{Li}_p(z)$ have been established in~\cite{OC}, alongside their global leading-order asymptotics for $\Re(p)>0$, our work focuses on the case $p=-\frac{1}{2}$. This specific value falls outside the validity of the asymptotic formula provided in~\cite{OC}. While we restrict our analysis to the negative real axis, we provide a significantly more precise description of the asymptotic behaviour in this regime.

\section{Computation of a fractional polylogarithm}\label{section2}

In an informative report~\cite{DW} on computation of polylogarithms, David C.\ Wood gave
formulae for the analytic continuation of a polylogarithm with $\operatorname{Li}_p(z)=\sum_{n=1}^\infty z^n/n^p$ for $|z|<1$. One obtains \eqref{Gsum} for $G(z)$ with $p=-\frac12$ as a particular case of~\cite[Eq.~(13.1)]{DW}.  We define $\log(z)$ on its first sheet, with $\Im(\log(z))\in(-\pi,\pi]$. 
Then the discontinuity of $\log(z)$ across its branch-cut on the negative
real $z$-axis causes no problem, since any integer multiple of $2\pi\im$ is absorbed by
shifting $n$ by an integer in the bilateral sum \eqref{Gsum}. For real $z\ge1$, 
the term with $n=0$ in \eqref{Gsum} creates a branch-cut.
We assume that $z\notin[1,\,\infty)$. Then \eqref{Gsum} has the complex-conjugation property
$\overline{G(z)}=G(\overline{z})$ that was requested by {}'t Hooft in~\cite[Sect. 3]{GH}. 

For the domain $\frac12<|z|<20<\e^{\pi}$, one can efficiently compute $G(z)$ using
\begin{gather}
G(z)=\frac{\sqrt{\pi}}{2\left(-\log(z)\right)^{3/2}}+\sum_{n=0}^\infty\zeta(-n-\tfrac12)\frac{(\log(z))^n}{n!},
\text{ for }\left|\log(z)\right|<2\pi,\label{Gzeta}\\[0.5ex]
\zeta(-n-\tfrac12)=-2\sin\!\left(\tfrac14(2n+1)\pi\right)\Gamma(n+\tfrac32)\,\zeta(n+\tfrac32)/(2\pi)^{n+3/2},\label{zeta}
\end{gather}
where \eqref{Gzeta} is a well-known result (see, e.g.,~\cite[Eq.~\href{http://dlmf.nist.gov/25.12.E12}{(25.12.12)}]{DLMF} or~\cite[Eq.~(9.3)]{DW}). The relation \eqref{zeta} follows from the analytic continuation of the Riemann zeta function, defined by $\zeta(s)=\sum_{n=1}^\infty n^{-s}$ for $\Re (s)>1$.

For $|z|>20$, it is preferable to use the inversion formula 
\begin{equation}
G(z)=\im\,G(1/z)+\frac{\im-1}{4\pi}\sum_{n=0}^\infty\left(n+\frac{\log(z)}{2\pi\im}\right)^{-3/2},
\text{ for }|z|\ge1\text{ and }\Im (z)\ge0,\label{Ginv}
\end{equation}
which is well known (see, for instance,~\cite[Eq.~\href{http://dlmf.nist.gov/25.12.E13}{(25.12.13)}]{DLMF} or~\cite[Eq.~(10.4)]{DW}). In particular, this relation yields the neat evaluation
\[
G(-1)=\frac{(1-2\sqrt2)\,\zeta(\frac32)}{4\pi}=-0.3801048126\ldots\,.
\]
With $|z|>20$ in \eqref{Ginv}, the first term on the right is quickly evaluated by \eqref{Gz} and in
the second term one encounters a complex Hurwitz zeta value, for which there is an efficient 
Euler--Maclaurin procedure~\cite{FJ}.

Most intriguingly, Wood gives an asymptotic expansion for $\operatorname{Li}_p(z)$
on the negative $z$-axis, with $z\ll-1$. Setting $z=-\e^{u}$, with large positive $u$,
and substituting $p=-\frac12$ in~\cite[Eq.~(11.1)]{DW}, one obtains the optimally truncated estimate
\begin{equation}
G(-\e^u)=-\frac{2}{\pi\sqrt{u}}\sum_{n=0}^{\lfloor u/2\rfloor}\eta(2n)\Gamma(2n+\tfrac12)u^{-2n}+\mathcal{O}(\e^{-u})\label{Gu} 
\end{equation}
where $\eta(s)=(1-2^{1-s})\zeta(s)$ and hence $\eta(0)=\frac12$. 

\section{A remarkable sequence of polynomials}\label{section3}

For real $u\ge0$, it follows from \eqref{Ginv} that
\begin{equation}
G(-\e^u)=\Re\!\left[\frac{\im-1}{4\pi}\sum_{n=0}^\infty\left(n+\frac12+\frac{u}{2\pi\im}\right)^{-3/2}\right]\label{Gh}
\end{equation}
is determined by a complex Hurwitz zeta value.
From \eqref{Gu}, it follows that for large $u$ one may express $G(-\e^u)$ in terms  of an optimally truncated sign-constant asymptotic expansion in $1/u^2$, together with an exponentially suppressed term.

For $u>0$, we define $S(u)$ by the Ansatz
\begin{equation}
G(-\e^u)=-\frac{2}{\pi\sqrt{u}}\left[\sum_{n=0}^{\lfloor u/2\rfloor}\eta(2n)\Gamma(2n+\tfrac12)u^{-2n}
+\sqrt{2\pi}\,\e^{-u}S(u)\right].\label{Ga}
\end{equation}
From numerical computation of \eqref{Gh},
we found that $S(u)\in(-0.7,\,0.4)$, for $u>0$, and that $S(u)=x-\frac23+\mathcal{O}(1/u)$, for large $u$,
with $x=u/2-\lfloor u/2\rfloor$. Moreover,
we found that this is the first term  of an asymptotic series, of the form 
\begin{equation}\label{Sexp}
S(u)=\sum_{k=0}^{\lfloor u\rfloor}\frac{P_k(x)}{u^k}+\mathcal{O}(\e^{-u}),\quad 
x=\tfrac{u}{2}-\left\lfloor\tfrac{u}{2}\right\rfloor\in[0,1),
\end{equation}
where $P_k(x)$ is a polynomial of degree $2k+1$ with rational coefficients. In the following, we first replicate the historical/empirical path pursued by the first author, which led to the discovery of these polynomials and highlighted their remarkable arithmetic structures. Subsequently, the second author provided the rigorous foundation, establishing the closed-form representation in Section~\ref{section4} which independently generates all these coefficients. Nonetheless, the first author's numerical approach using high-precision arithmetic remains crucial for validating the boundaries of our asymptotic expansion and guiding the formulation of Theorem~\ref{thm1}. Consequently, the rest of Section~\ref{section3} can be safely bypassed by readers who wish to focus exclusively on the rigorous development presented from Section~\ref{section4} onwards.

The sequence of polynomials begins with
\begin{align*}
P_0(x)&=x -\tfrac23, \\[0.5ex]
P_1(x)&=\tfrac23x^3 - x^2 + \tfrac{7}{24}x + \tfrac{47}{2160}, \\[0.5ex]
P_2(x)&=\tfrac25x^5 - \tfrac23x^4 - \tfrac{1}{36}x^3 + \tfrac13x^2 - \tfrac{73}{1920}x - \tfrac{433}{24192}, \\[0.5ex]
P_3(x)&=\tfrac{4}{21}x^7 - \tfrac29x^6 - \tfrac{5}{12}x^5 + \tfrac{31}{72}x^4 
+ \tfrac{433}{1728}x^3 - \tfrac{223}{1152}x^2 - \tfrac{106619}{2903040}x + \tfrac{28583}{2488320}, 
\end{align*}
and has been developed up to $k=166$. 
The denominators of the coefficients of $P_k(x)$ involve no prime greater than $2k+3$.

The differences $\Delta_k(x)=P_k(x+1)-P_k(x)$ are determined by the asymptotic series
\[
\frac{\e^u}{\sqrt{2\pi}}\,\frac{\Gamma(u-2x+\tfrac12)}{u^{u-2x}} \sim \sum_{k=0}^\infty\frac{\Delta_k(x)}{u^k}
=1+\frac{2x^2 - \tfrac{1}{24}}{u}+\mathcal{O}\!\left(\frac{1}{u^2}\right),
\]
since the transformation  $x\to x+1$ would correspond to the instruction to omit the
last term of the summation in \eqref{Ga}, at $n=\frac12u-x$. Taking a logarithm, we obtain
\begin{equation}
\log\left[\sum_{k=0}^\infty\frac{\Delta_k(x)}{u^k}\right] \sim
\sum_{n = 2}^\infty  \frac{B_n \!\left( 2x + \tfrac{1}{2} \right)}{n(n - 1)}\frac{1}{u^{n - 1} },\label{dlog}
\end{equation}
where $B_n\!\left( 2x + \frac{1}{2} \right)$ is a Bernoulli polynomial (cf.~\cite[Eqs.~\href{http://dlmf.nist.gov/5.11.E8}{(5.11.8)} and \href{http://dlmf.nist.gov/24.4.E3}{(24.4.3)}]{DLMF}).

The finite difference equation $P_k(x+1)=P_k(x)+\Delta_k(x)$ determines the polynomial $P_k(x)$
modulo its constant term, $P_k(0)$. To determine $P_k(0)$, the first author resorted to experiment, using 
{\tt zetahurwitz} in {\tt Pari/GP}~\cite{Pari} 
to compute instances of $G(-\e^u)$ at 600  even integers, $u_n\in[6002,\,7200]$,
working at 6000-digit precision. This took 30 minutes on  a single core.
Then \eqref{Ga} gives sufficient precision to determine, iteratively, 167 rational values of $P_k(0)$
from 600 expansions $S(u_n)\approx\sum_{k=0}^{166}P_k(0)/u_n^k$.
The iterative process relies on control of the denominator $D_k$ of $P_k(0)=N_k/D_k$.
We found that $D_k/D_{k-1}$ is a relatively small rational number,
involving no prime greater than $2k+3$. For example 
$D_{166}/D_{165}=2^3\cdot3^3\cdot5\cdot11\cdot113=1342440$, while $N_{166}$
is a 780-digit integer, obtained from numerical data with an absolute error less than $10^{-27}$
and hence with good confidence. 

Combining this empirical data for $P_k(0)$ with difference polynomials $\Delta_k(x)$, obtained 
by rational linear algebra from \eqref{dlog}, we determined $P_k(x)$ exactly for $k\in[0,166]$. 
Studying these results we found that $P_k(x)$ at large $k$ is very well approximated 
by a sinusoid whose amplitude grows with $k$ exponentially faster than that for $\Delta_k(x)$.
Moreover, the phase of this sinusoid for $P_k(x)$ increases linearly with $k$.

From \eqref{dlog}, one sees that  $g_k=\Delta_k(-\frac14)$ has a generating function~\cite{WB,GN}
\begin{equation}\label{ggenf}
\sum_{k=0}^\infty g_ky^k=\exp\!\left(\sum_{n=1}^\infty\frac{B_{2n}\,y^{2n-1}}{2n(2n-1)}\right)=
1+\tfrac{1}{12}y+\tfrac{1}{288}y^2-\tfrac{139}{51840}y^3-\tfrac{571}{2488320}y^4
+\mathcal{O}(y^5).
\end{equation}
At large $k$, one has $kg_k=\mathcal{O}(k!/(2\pi)^k)$. The detailed behaviour depends on the parity of $k$, with resurgent asymptotic expansions given by~\cite{WB}
\begin{equation}
g_{2m}\sim-2\sum_{n=0}^{\lfloor m/2\rfloor}g_{2n+1}\frac{\Gamma(2m-2n-1)}{(2\pi\im)^{2m-2n}},\quad
g_{2m-1}\sim-2\sum_{n=0}^{\lfloor m/2\rfloor}g_{2n}\frac{\Gamma(2m-2n-1)}{(2\pi\im)^{2m-2n}}\label{gres}
\end{equation}
optimally truncated at $n=\lfloor m/2\rfloor$. We found that 
\begin{equation}
\Delta_k(x) = 
\frac{2 \, \Gamma(k)}{\left(2\pi\right)^{k+1}}
\left(\sin\!\left(4\pi x - \tfrac{\pi}{2}k\right)
  + \mathcal{O}\!\left(\frac{1}{k}\right)
\right),\label{Dkasy}
\end{equation}
for large $k$ and bounded real $x$ (see Appendix~\ref{Dkasyproof}). At $x=-\frac14$, this accords with the leading terms in \eqref{gres}.

Solving the difference equation $P_k(x+1)=P_k(x)+\Delta_k(x)$, with
the boundary value $P_k(0)$ determined from fits to Hurwitz zeta values, we found a very different 
sinusoidal pattern for $P_k(x)$ at large $k$.

\section{Asymptotic constants}\label{section4}
 
For $x\in[0,1)$ and large $k$, there are positive constants $(C,R)$ 
such that
\begin{equation}
P_k(x)=\frac{1}{\sqrt{2\pi}}\left(\sin((2k+1)C-2\pi x)
+\mathcal{O}\!\left(\frac{1}{\sqrt{k}\left(2\pi R^2\right)^k}\right)\right)R^{2k+1}\Gamma\!\left(k+\tfrac12\right).\label{Pkasy}
\end{equation}
From exact results for $P_k(x)$, with $k\le100$, we obtained the approximate values
\[
C\approx1.0688539158679530121571,\quad R\approx0.5181839789815558726739.
\]
The correction to the sinusoid in \eqref{Pkasy} is  suppressed by a factor 
$k^{-1/2}\exp(-Dk)$, with $D=\log(2\pi R^2)>0.523$, while $C$ determines the rate at which the phase of $P_k(x)$ 
increases with $k$. The frequency of the sinusoid in the accurate formula for $P_k(x)$ at large $k$ is half the frequency of  the rough approximation of $\Delta_k(x)$ in \eqref{Dkasy}. 

Remarkably, one does not need to evaluate more Hurwitz zeta values to improve the estimates for $C$ and $R$,
since the derivatives $P_k^\prime(0)$ and $P_k^{\prime\prime}(0)$ suffice for this purpose.
These are determined by $\Delta_k(x)=P_k(x+1)-P_k(x)$, using rational linear algebra. 
Performing this algebraic task up to $k=450$, the first author obtained 100 good digits of $C$ and $R$.
Then the second author identified these constants as follows.

\begin{theorem}\label{thm1}
The polynomials $P_k(x)$ obey \eqref{Pkasy}
as $k\to+\infty$, uniformly for bounded real values of $x$, with positive constants $C$ and $R$ satisfying
\begin{equation}
R\exp(\im\,C)=\left(\frac{-1}{2+\pi\im}\right)^{1/2}.\label{Th1}
\end{equation}
In particular,
\[
C = \tfrac{\pi}{2} - \tfrac{1}{2}\arctan\!\left(\tfrac{\pi}{2}\right),\quad  R = \left(4 + \pi^{2}\right)^{-1/4}.
\]
\end{theorem}

The proof of \eqref{Pkasy} and Theorem~\ref{thm1} relies on the following representation of the polynomials $P_k(x)$:
\begin{gather}\label{Pkint}
\begin{split}
P_k(x) &= \frac{1}{\sqrt{2\pi}} \frac{\Gamma\!\left(k + \frac{1}{2}\right)}{2\pi \im} \oint_{|t|=1} 
\frac{\e^{(5/2- 2x) t}}{1 - \e^{2t}} 
\frac{\d t}{\left(\e^t - t - 1\right)^{k + 1/2}}
 \\ &= \frac{1}{2^k(2k+1) k!}\left[ \frac{\d^{2k+1} }{\d t^{2k+1} }\left( \frac{t\,\e^{(5/2 - 2x)t} }{1 - \e^{2t} }\left( \frac{1}{2}\frac{t^2 }{\e^t  - t - 1} \right)^{k + 1/2}  \right) \right]_{t = 0} ,
\end{split}
\end{gather}
which will be derived in the next section. From the second representation, it is seen that the $P_k(x)$ are polynomials in $x$ of degree $2k+1$ with rational coefficients whose denominators involve no prime greater than $2k+3$. The interested reader may verify this arithmetic property by substituting the standard Maclaurin series for $(1 - \e^{2t})/t$ and $2(\e^t - t - 1)/t^2$ into the representation, expanding the reciprocal and the negative $(k+1/2)^{\mathrm{th}}$ power via the geometric and binomial series up to order $2k+1$, and tracking the maximum primes in the resulting factorials. The bound $2k+3$ then follows naturally from the index limits of these intermediate expansions and the external factor.

\section{Proofs}\label{section5}

The main goal of this section is to prove Theorem~\ref{thm1}. We begin by establishing the general expansion \eqref{expansion} for $G(-\e^u)$, expressed in terms of terminant functions (see Appendix~\ref{terminant} for their definition). This expansion leads to the representation \eqref{Sformula} for $S(u)$, and subsequently to the integral formula \eqref{Pkint} for the polynomials $P_k(x)$. The remainder of the section is then devoted to establishing the asymptotic formula \eqref{Pkasy}.

By analytically continuing~\cite[Eq.~(3.3)]{TMG} to $p\le -1$ and substituting $p=-\frac{3}{2}$, we can write
\[
- \frac{\pi \sqrt{u}}{2} G(-\e^u) 
= \frac{\sqrt{\pi}}{2} 
+ \frac{\sqrt{\pi}}{4} \sum_{k=1}^\infty (-1)^k U\!\left(1,\tfrac{1}{2},ku\right) 
+ \frac{\sqrt{\pi}}{2} \sum_{k=1}^\infty (-1)^k M\!\left(1,\tfrac{1}{2},-ku\right),
\]
where $M$ and $U$ are the confluent hypergeometric functions~\cite[\href{http://dlmf.nist.gov/13.2}{\S13.2}]{DLMF}. Using~\cite[\href{http://dlmf.nist.gov/13.6.ii}{\S13.6(ii)}]{DLMF}, these functions can be expressed in terms of incomplete gamma functions:
\[
U\!\left(1, \tfrac{1}{2}, ku \right) = \sqrt{ku} \,\e^{ku} \Gamma\!\left(-\tfrac{1}{2}, ku\right),
\]
and
\[
M\!\left(1, \tfrac{1}{2}, -ku \right) 
= \operatorname{Re}\!\left[ -\tfrac{\im}{2} \sqrt{ku} \,\e^{-ku} \gamma\!\left(-\tfrac{1}{2}, k u \e^{\pi \im} \right) \right] 
= -\tfrac{1}{2} \sqrt{ku} \,\e^{-ku} \operatorname{Im} \Gamma\!\left(-\tfrac{1}{2}, k u \e^{\pi \im} \right),
\]
for $u>0$ and $k \in \mathbb{Z}_{>0}$. An application of~\cite[Eq.~\href{http://dlmf.nist.gov/8.8.E10}{(8.2.9)}]{DLMF} then gives the following expansions:
\[
U\!\left(1,\tfrac{1}{2},ku\right) 
= \frac{2}{\sqrt{\pi}} \sum_{n=1}^{2N_k} (-1)^{n+1} \frac{1}{k^n} 
   \Gamma\!\left(n+\tfrac{1}{2}\right) u^{-n} 
   - 4 \sqrt{\pi ku}\, \e^{ku} T_{2N_k+3/2}(ku),
\]
and
\[
M\!\left(1,\tfrac{1}{2},-ku\right) 
= - \frac{1}{\sqrt{\pi}} \sum_{n=1}^{2N_k} \frac{1}{k^n} 
   \Gamma\!\left(n+\tfrac{1}{2}\right) u^{-n} 
   + 2 \sqrt{\pi ku} \,\e^{-ku} \Im T_{2N_k+3/2}\!\left(ku\e^{\pi \im}\right),
\]
where the $N_k$ are arbitrary non-negative integers and $T_p(z)$ is defined by \eqref{Tpdef}. Accordingly, we arrive at the exact expansion
\begin{gather}\label{expansion}
\begin{split}
&- \frac{\pi \sqrt{u}}{2} G(-\e^u) 
= \frac{\sqrt \pi}{2}  + \sum_{k = 1}^\infty \sum_{n = 1}^{N_k } \frac{( - 1)^{k + 1} }{k^{2n}}\Gamma\!\left(2n+\tfrac{1}{2}\right)u^{-2n}
\\ &\quad - \pi \sqrt{u} \sum_{k=1}^\infty (-1)^k \sqrt{k} \,\e^{ku} T_{2N_k+3/2}(ku)
 + \pi \sqrt{u} \sum_{k=1}^\infty (-1)^k \sqrt{k} \,\e^{-ku} 
   \Im T_{2N_k+3/2}\!\left(ku\e^{\pi \im}\right),
\end{split}
\end{gather}
which holds for $u>0$ and for arbitrary non-negative integers $N_k$. Truncating each $n$-series near its least term, given approximately by $N_k \approx ku$, yields an exponentially improved expansion, closely analogous to the one known for the logarithm of the gamma function~\cite{MVB}.

If we set $N_k = N$ for all $k$, then the order of summation in \eqref{expansion} may be interchanged, yielding
\begin{align*}
- \frac{\pi \sqrt{u}}{2} G(-\e^u) 
= \sum_{n=0}^N \eta(2n) \Gamma\!\left(2n+\tfrac{1}{2}\right) u^{-2n} 
& - \pi \sqrt{u} \sum_{k=1}^\infty (-1)^k \sqrt{k} \,\e^{ku} T_{2N+3/2}(ku)
\\ & + \pi \sqrt{u} \sum_{k=1}^\infty (-1)^k \sqrt{k}\, \e^{-ku} 
   \Im T_{2N+3/2}\!\left(ku\e^{\pi \im}\right).
\end{align*}
Applying the asymptotic expansions \eqref{asymp1} and \eqref{asymp2} for $T_p(z)$, we obtain
\begin{align*}
- \frac{\pi \sqrt{u}}{2} G(-\e^u) 
& = \sum_{n=0}^N \eta(2n) \Gamma\! \left(2n+\tfrac{1}{2}\right) u^{-2n} 
\\ &\quad + \sqrt{2\pi}\, \e^{-u} 
\left[ \sqrt{\tfrac{\pi u}{2}} \e^{2u} T_{2N+3/2}(u) 
- \sqrt{\tfrac{\pi u}{2}} \Im T_{2N+3/2} \!\left(u\e^{\pi \im}\right) \right] 
+ \mathcal{O}(\e^{-2u}),
\end{align*}
as $u \to +\infty$ with $N \approx u/2$. Consequently, by \eqref{Ga}, we obtain
\begin{equation}\label{Sformula}
S(u)= \sqrt{\tfrac{\pi u}{2}} \e^{2u} T_{2\lfloor u/2\rfloor+3/2}(u) 
- \sqrt{\tfrac{\pi u}{2}} \Im T_{2\lfloor u/2\rfloor+3/2} \!\left(u\e^{\pi \im}\right) + \mathcal{O}(\e^{-u}),
\end{equation}
as $u \to +\infty$. Using the asymptotic expansions \eqref{asymp1} and \eqref{asymp2}, it follows that
\begin{equation}\label{asymp3}
\sqrt{\tfrac{\pi u}{2}} \e^{2u} T_{2\lfloor u/2\rfloor+3/2}(u) 
- \sqrt{\tfrac{\pi u}{2}} \Im T_{2\lfloor u/2\rfloor+3/2} \!\left(u\e^{\pi \im}\right)\sim \sum_{k=0}^\infty \frac{P_k(x)}{u^k},
\end{equation}
as $u\to+\infty$, where $x = u/2 - \lfloor u/2 \rfloor$. The coefficients are given explicitly by
\begin{align*}
P_k(x) &= - \frac{1}{2\sqrt{2\pi}} 
\left( A_k\!\left(\tfrac{3}{2}-2x\right) - B_k\!\left(\tfrac{3}{2}-2x\right) \right)
 \\ &= \frac{1}{\sqrt{2\pi}} \frac{\Gamma\!\left(k + \frac{1}{2}\right)}{2\pi \im} \oint_{|t|=1} 
\frac{\e^{(5/2- 2x) t}}{1 - \e^{2t}} 
\frac{\d t}{\left(\e^t - t - 1\right)^{k + 1/2}}
 \\ &= \frac{1}{2^k(2k+1) k!}\left[ \frac{\d^{2k+1} }{\d t^{2k+1} }\left( \frac{t\,\e^{(5/2 - 2x)t} }{1 - \e^{2t} }\left( \frac{1}{2}\frac{t^2 }{\e^t  - t - 1} \right)^{k + 1/2}  \right) \right]_{t = 0} ,
\end{align*}
valid for all $k \ge 0$. By the uniqueness of the coefficients in an asymptotic expansion, these coefficients coincide with the polynomials appearing in \eqref{Sexp}.

We now justify the error term in \eqref{Sexp}, assuming that Theorem~\ref{thm1} has been established. From \eqref{Pkasy}, the optimal truncation of the asymptotic series \eqref{asymp3} occurs near $k \approx \lfloor u / R^2 \rfloor$. The error introduced by truncating the expansion at or around its least term is of the same order of magnitude as the first omitted term, multiplied by an additional factor of $\sqrt{u}$ arising from the Stokes phenomenon (see, e.g.,~\cite[Sect. 8]{WB1} or~\cite[Sect. 5]{ODO}), since the positive real axis is a Stokes line for the left-hand side of \eqref{asymp3}. Accordingly,
\begin{align*}
S(u) 
&= \sum_{k=0}^{\lfloor u / R^2 \rfloor} \frac{P_k(x)}{u^k} 
   + \mathcal{O}\!\left(\sqrt{u}\,\e^{-u / R^2}\right) + \mathcal{O}(\e^{-u}) \\
&= \sum_{k=0}^{\lfloor u \rfloor} \frac{P_k(x)}{u^k} 
   + \sum_{k=\lfloor u \rfloor + 1}^{\lfloor u / R^2 \rfloor} \frac{P_k(x)}{u^k} 
   + \mathcal{O}\!\left(\sqrt{u}\,\e^{-u / R^2}\right) + \mathcal{O}(\e^{-u})= \sum_{k=0}^{\lfloor u \rfloor} \frac{P_k(x)}{u^k} 
   + \mathcal{O}(\e^{-u}),
\end{align*}
as $u \to +\infty$, thereby establishing \eqref{Sexp}. We note that it actually suffices to sum up to $\lfloor \lambda u \rfloor$ in \eqref{Sexp}, where $\lambda = -1/W_{-1}\!\left(-R^2/\e\right) \approx 0.2782$, and $W_{-1}$ denotes the non-principal real branch of the Lambert $W$-function~\cite[\href{http://dlmf.nist.gov/4.13}{\S4.13}]{DLMF}.

\begin{figure}[!ht]
  \centering
  \includegraphics[width=0.4\textwidth]{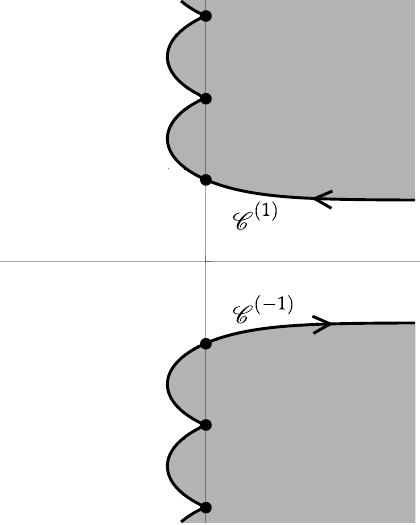}
  \caption{Illustration of the domain $D$ (unshaded) in the $t$-plane, bounded by the steepest descent contours $\mathscr{C}^{(\pm1)}$. The dots indicate the critical points $\pm2\pi \im, \pm4\pi \im,\pm6\pi\im,\ldots$.}
  \label{fig1}
\end{figure}

We conclude this section with the proof of Theorem~\ref{thm1}. We first examine the analytic properties of the function $t\mapsto \sqrt{\e^t-t-1}$ appearing in the integral representation \eqref{Pkint}. Let $\mathscr{P}(\theta)$ denote the steepest descent path emanating from the origin in the $t$-plane along which
\[
\arg \!\left( \e^{\im\theta } \left( \e^t  - t - 1 \right) \right) = 0,
\]
and which continuously deforms, as $\theta$ varies, from the positive real axis corresponding to $\theta=0$. As $\theta$ increases from $0$ to $4\pi$, these paths sweep out a domain in the $t$-plane, denoted by $D$. The contour $\mathscr{P}(\theta)$ changes discontinuously at $\theta=\frac{5\pi}{2}$ and $\theta=\frac{3\pi}{2}$, where it connects to the critical points $2\pi \im, 4\pi \im, 6\pi \im,\ldots $ and $-2\pi \im, -4\pi \im, -6\pi \im,\ldots $ of $t\mapsto \e^t-t-1$, respectively. The boundary of the domain $D$ consists of two infinite contours, $\mathscr{C}^{(\pm1)}$. Each of these is itself a steepest descent path passing through the corresponding critical point $\pm2\pi \im$, along which
\[
\arg \!\left( \e^{ \pm \frac{\pi }{2}\im} \left( \e^t  - t-1 \pm 2\pi\im \right) \right) = 0 \bmod 2\pi
\]
for $t\in \mathscr{C}^{(\pm1)}$. The domain $D$ and the contours $\mathscr{C}^{(\pm1)}$ are illustrated in Figure~\ref{fig1}. We fix the phase of $\sqrt{\e^t-t-1}$ along the contour $\mathscr{P}(\theta)$ by defining
\begin{equation}\label{sqrtdef}
\arg\! \left( \sqrt {\e^t  - t - 1} \right) =  - \tfrac{\theta }{2}.
\end{equation}
Under this convention, the function $t\mapsto \sqrt{\e^t-t-1}$ is analytic throughout the domain $D$ and continuous on its boundary.

We begin by deforming the loop contour in \eqref{Pkint} into a circle of radius $5$, which lies entirely within the domain $D$. During this deformation, the contour crosses the poles at $t = \pm \pi \im$, contributing the corresponding residues. We thus obtain
\begin{align*}
P_k(x) &= \frac{\im}{2\sqrt{2\pi}} \Gamma\!\left(k + \tfrac{1}{2}\right)\frac{\e^{- 2\pi x \im} }{\left( - 2 - \pi \im\right)^{k + 1/2} } - \frac{\im}{2\sqrt{2\pi}} \Gamma\!\left(k + \tfrac{1}{2}\right)\frac{\e^{ 2\pi x \im} }{\left( - 2 + \pi \im\right)^{k + 1/2} } \\ &\quad + \frac{1}{\sqrt{2\pi}} \frac{\Gamma\!\left(k + \frac{1}{2}\right)}{2\pi \im} \oint_{|t|=5} 
\frac{\e^{(5/2- 2x) t}}{1 - \e^{2t}} 
\frac{\d t}{\left(\e^t - t - 1\right)^{k + 1/2}}. 
\end{align*}
The first instance at which $\mathscr{P}(\theta)$ reaches $\pi \im$ and $-\pi \im$ occurs at $\theta=\pi+\arctan\! \left( \tfrac{\pi }{2} \right)$ and $\theta=3\pi-\arctan\! \left( \tfrac{\pi }{2} \right)$, respectively. Hence, by \eqref{sqrtdef}, we have
\[
\arg \!\left(\sqrt{ - 2 \mp \pi \im}\right) = \arg \!\left( \sqrt{\e^{\pm\pi \im}  \mp \pi \im - 1} \right) = \pm\tfrac{\pi}{2}  \pm \tfrac{1}{2}\arctan\! \left( \tfrac{\pi }{2} \right) \bmod 2\pi.
\]
Accordingly, we can write
\begin{align*}
P_k(x) &=\frac{1}{\sqrt{2\pi}} \sin((2k+1)C-2\pi x)
 R^{2k+1}\Gamma\!\left(k+\tfrac12\right)
\\ &\quad + \frac{1}{\sqrt{2\pi}} \frac{\Gamma\!\left(k + \frac{1}{2}\right)}{2\pi \im} \oint_{|t|=5} 
\frac{\e^{(5/2- 2x) t}}{1 - \e^{2t}} 
\frac{\d t}{\left(\e^t - t - 1\right)^{k + 1/2}} .
\end{align*}
Therefore, it remains to establish that the contour integral is of order $\mathcal{O}(k^{-1/2}/(2\pi)^k)$ for large $k$. To do this, we rewrite the integral in an alternative form. For $1 \le \beta < 2$, define
\begin{equation}\label{Fdef}
F_k(x,\beta)= \frac{1}{\sqrt{2\pi}} \frac{\Gamma\!\left(k + \frac{1}{2}\right)}{2\pi \im} \oint_{|t|=5} 
\frac{\e^{(5/2- 2x) t}}{\beta^2 - \e^{2t}} 
\frac{\d t}{\left(\e^t - t - 1\right)^{k + 1/2}}.
\end{equation}
With this notation, we have
\begin{equation}\label{Pleading}
P_k(x) =\frac{1}{\sqrt{2\pi}} \sin((2k+1)C-2\pi x)R^{2k+1}\Gamma\!\left(k+\tfrac12\right)+F_k(x,1).
\end{equation}
Assume that $\frac{3-k}{2} < x < \frac{5}{4}$, $k\ge1$, and $1 < \beta < 2$. Under these conditions, we first derive an alternative form of $F_k(x,\beta)$. An alternative form of $F_k(x,1)$ is then derived by passing to the limit $\beta \to 1^+$.

We begin by deforming the circle of integration in \eqref{Fdef} onto the boundary of the domain $D$. The contributions from the arcs at infinity vanish in view of our assumptions that $\frac{3-k}{2} < x < \frac{5}{4}$ and $k\ge1$. This leads to
\[
F_k(x,\beta)= \frac{1}{\sqrt{2\pi}}  \sum_{m=\pm1}\frac{\Gamma\!\left(k + \frac{1}{2}\right)}{2\pi \im}\int_{\mathscr{C}^{(m)}} 
\frac{\e^{(5/2- 2x) t}}{\beta^2 - \e^{2t}} 
\frac{\d t}{\left(\e^t - t - 1\right)^{k + 1/2}}.
\]
Note that, since $\beta > 1$, all poles of the integrand lie outside the region through which the contour was deformed. Next, introducing the Euler integral representation of the gamma function, we obtain
\begin{equation}\label{F2}
F_k(x,\beta)= \frac{1}{\sqrt{2\pi}}  \sum_{m=\pm1}\int_{0}^{+\infty}\e^{-s}s^{k-1/2}\frac{1}{2\pi \im}\int_{\mathscr{C}^{(m)}} 
\frac{\e^{(5/2- 2x) t}}{\beta^2 - \e^{2t}} 
\frac{\d t}{\left(\e^t - t - 1\right)^{k + 1/2}}\id s.
\end{equation}
Following the procedure in~\cite[p. 615]{WB}, we change variables from $s$ and $t$ to $v$
and $t$ via
\begin{equation}\label{change}
s = \pm \im v\left(\e^t  - t - 1\right).
\end{equation}
When applying \eqref{change} to the integrals in \eqref{F2}, the square roots in the integrands should be interpreted as
\[
\sqrt { \pm  \im}  = \e^{ \mp \frac{3\pi}{4}\im}. 
\]
With these substitutions, the representation \eqref{F2} becomes
\begin{gather}\label{F3}
\begin{split}
F_k(x,\beta)&= \frac{\e^{-\frac{3\pi \im}{4}}}{\sqrt{2\pi}}\im^k
\int_{0}^{+\infty} \e^{-2\pi v} v^{k - 1/2}\frac{1}{2\pi \im}
  \int_{\mathscr{C}^{(1)}}
  \e^{-\im v(\e^t - t - 1 + 2\pi \im)}
  \frac{\e^{(5/2 - 2x)t}}{\beta^2 - \e^{2t}} \id t
 \id v \\
&\quad
+ \frac{\e^{\frac{3\pi \im}{4}}}{\sqrt{2\pi}} (-\im)^k
\int_{0}^{+\infty} \e^{-2\pi v} v^{k - 1/2}\frac{1}{2\pi \im}
  \int_{\mathscr{C}^{(-1)}}
  \e^{\im v(\e^t - t - 1 - 2\pi \im)}
  \frac{\e^{(5/2 - 2x)t}}{\beta^2 - \e^{2t}} \id t
 \id v.
\end{split}
\end{gather}
We define the contour $\overline{\mathscr{C}}^{(1)}$ to be congruent to $\mathscr{C}^{(1)}$ but shifted downward in the complex plane by $2 \pi \im$, and similarly, $\overline{\mathscr{C}}^{(-1)}$ to be congruent to $\mathscr{C}^{(-1)}$ but shifted upward by $2 \pi \im$. With these contours, the representation \eqref{F3} become
\begin{gather}\label{F4}
\begin{split}
F_k(x,\beta)&=\frac{\e^{-\frac{3\pi \im}{4}} \e^{2\pi \im(5/2 - 2x)}}{\sqrt{2\pi}} 
\im^k
\int_{0}^{+\infty} \e^{-2\pi v} v^{k - 1/2}\frac{1}{2\pi \im}
\int_{\overline{\mathscr{C}}^{(1)}}
\e^{-\im v(\e^t - t - 1)}
\frac{\e^{(5/2 - 2x)t}}{\beta^2 - \e^{2t}} \id t
\id v \\ 
&\quad
+ \frac{\e^{\frac{3\pi \im}{4}} \e^{-2\pi \im(5/2 - 2x)}}{\sqrt{2\pi}} 
(-\im)^k
\int_{0}^{+\infty} \e^{-2\pi v} v^{k - 1/2}\frac{1}{2\pi \im}
\int_{\overline{\mathscr{C}}^{(-1)}}
\e^{\im v(\e^t - t - 1)}
\frac{\e^{(5/2 - 2x)t}}{\beta^2 - \e^{2t}} \id t
\id v.
\end{split}
\end{gather}
Next, we rewrite the contour integrals in terms of the terminant function $T_p(z)$. To see this, note that from \eqref{Tintegrals} and \eqref{Tpoleint} we have
\[
\frac{\beta^{z+\alpha-2} \e^z}{2} \left[
    \e^{-\pi \im (z+\alpha)} \e^{\beta z} T_{z+\alpha}(\beta z)
    + \e^{-\beta z} T_{z+\alpha}\!\left( \beta z \e^{-\pi \im} \right)
\right]
= \frac{1}{2\pi \im} \cupint_{-\infty}^{+\infty} 
    \e^{-z(\e^t - t - 1)} \frac{\e^{\alpha t}}{\beta^2 - \e^{2t}} \d t,
\]
provided that $\beta\ge 1$, $\Re(z+\alpha)>0$, and $\left|\arg z\right|<\frac{\pi}{2}$. The integration contour is chosen to pass below the singularity at $t=\log(\beta)$. This representation can be continued analytically in $z$ by deforming the contour to follow a path of steepest descent while avoiding the singularity at $t=\log(\beta)$. In this way, we obtain
\begin{gather}\label{contint1}
\begin{split}
&\frac{1}{2\pi \im}
\int_{\overline{\mathscr{C}}^{(1)}}
\e^{-\im v(\e^t - t - 1)}
\frac{\e^{(5/2 - 2x)t}}{\beta^2 - \e^{2t}} \id t \\ & \quad =\frac{\beta^{\im v + 1/2 - 2x} \e^{\im v}}{2}
\left[
  \e^{-\pi \im (\im v + 5/2 - 2x)} \e^{\beta \im v}
  T_{\im v + 5/2 - 2x}\!\left( \beta v \e^{\frac{\pi }{2}\im} \right)
  + \e^{-\beta \im v}
  T_{\im v + 5/2 - 2x}\!\left( \beta v \e^{-\frac{\pi }{2}\im} \right)
\right].
\end{split}
\end{gather}
Taking complex conjugates and applying \eqref{Tconj}, we obtain the corresponding expression for the other contour integral:
\begin{gather}\label{contint2}
\begin{split}
&\frac{1}{2\pi \im}
\int_{\overline{\mathscr{C}}^{(-1)}}
\e^{\im v(\e^t - t - 1)}
\frac{\e^{(5/2 - 2x)t}}{\beta^2 - \e^{2t}} \id t =-\frac{\beta^{-\im v + 1/2 - 2x} \e^{-\im v}}{2}\\ & \quad \times
\left[
  \e^{-\pi \im (-\im v + 5/2 - 2x)} \e^{-\beta \im v}
  T_{-\im v + 5/2 - 2x}\!\left( \beta v \e^{-\frac{\pi}{2}\im} \right)
  + \e^{-2\pi \im (-\im v + 5/2- 2x)} \e^{\beta \im v}
  T_{-\im v + 5/2 - 2x}\!\left( \beta v \e^{\frac{\pi}{2}\im} \right)
\right].
\end{split}
\end{gather}
Using the integral representation \eqref{Tintegral}, the reader can verify that the right-hand sides of \eqref{contint1} and \eqref{contint2} are both
\[
\mathcal{O}\big( v^{ - 5/2 + 2x} \e^{\frac{\pi }{2}v}  \big),
\]
uniformly for $1 < \beta < 2$, real $x < \tfrac{5}{4}$ bounded, and all $v > 0$. Thus, we may substitute the right-hand sides of \eqref{contint1} and \eqref{contint2} into that of \eqref{F4}, and then pass to the limit $\beta \to 1^+$ by applying the dominated convergence theorem. This yields
\begin{gather}\label{F1expr}
\begin{split}
& F_k(x,1)=\frac{\e^{-\frac{3\pi \im}{4}} \e^{2\pi \im(5/2 - 2x)}}{2\sqrt{2\pi}} 
\im^k
\int_{0}^{+\infty} \e^{-2\pi v} v^{k - 1/2} \\ & \qquad\times 
\left(
  \e^{-\pi \im (\im v + 5/2 - 2x)} \e^{2 \im v}
  T_{\im v + 5/2 - 2x}\!\left( v \e^{\frac{\pi}{2}\im} \right)
  + T_{\im v + 5/2 - 2x}\!\left( v \e^{-\frac{\pi }{2}\im} \right)
\right)
\d v \\ 
&\quad
- \frac{\e^{\frac{3\pi \im}{4}} \e^{-2\pi \im(5/2 - 2x)}}{2\sqrt{2\pi}} 
(-\im)^k
\int_{0}^{+\infty} \e^{-2\pi v} v^{k -1/2} \\ & \qquad\times 
\left( \e^{ - \pi \im( - \im v + 5/2 - 2x)} \e^{ - 2\im v} T_{ - \im v + 5/2 - 2x}\! \left( v\e^{ - \frac{\pi }{2}\im}  \right) + \e^{ - 2\pi \im( - \im v + 5/2 - 2x)} T_{ - \im v + 5/2 - 2x} \!\left( v\e^{\frac{\pi }{2}\im} \right) \right)
\d v,
\end{split}
\end{gather}
provided that $\tfrac{3-k}{2} < x < \tfrac{5}{4}$. For large $k$, the dominant contribution to the integrals arises from large values of $v$. Taking into account the remark following \eqref{Acoeff} and using \eqref{asymp1}, we find that
\[
\e^{ -\pi \im( \pm \im v + 5/2 - 2x)} \e^{ \pm 2\im v} T_{ \pm \im v + 5/2 - 2x} \!\left( v\e^{ \pm \frac{\pi }{2}\im}  \right) = \mathcal{O}\big( v^{ - 1/2} \big)
\]
as $v\to+\infty$, uniformly for bounded real $x<\frac{5}{4}$. Using the integral representation \eqref{Tintegral}, one readily verifies that
\begin{equation}\label{Texpsmall}
T_{\im v + 5/2 - 2x} \!\left( v\e^{ - \frac{\pi }{2}\im} \right)= \mathcal{O}\big( \e^{ - \frac{3\pi}{2}v}  \big),\quad \e^{ - 2\pi \im( - \im v + 5/2 - 2x)} T_{ - \im v + 5/2 - 2x} \!\left( v\e^{\frac{\pi }{2}\im} \right) = \mathcal{O}\big( \e^{ - \frac{3\pi}{2}v}  \big)
\end{equation}
as $v\to+\infty$, uniformly for bounded real $x<\frac{5}{4}$. Accordingly, we conclude that
\[
F_k (x,1) = \mathcal{O}\!\left( \frac{\Gamma (k)}{\left(2\pi \right)^k } \right),
\]
again uniformly for bounded real $x<\frac{5}{4}$. By substituting into \eqref{Pleading} and employing $\Gamma\!\left(k+\frac{1}{2}\right)\sim \sqrt{k}\,\Gamma(k)$~\cite[Eq.~\href{http://dlmf.nist.gov/5.11.E12}{(5.11.12)}]{DLMF}, Theorem~\ref{thm1} follows under the added restriction $x<\frac{5}{4}$. To remove this restriction, we can proceed by induction using the difference equation $P_k (x + 1) = P_k (x) + \Delta _k (x)$ together with the rough asymptotics \eqref{Dkasy} for $ \Delta _k (x)$. This completes the proof of Theorem~\ref{thm1}.

\section{Discussion}\label{section6}

In this paper, we have analysed the analytic continuation and asymptotic behaviour of the series $G(z)=\sum_{n=1}^\infty\sqrt{n}\,z^n$ originally proposed by 't Hooft in the context of the quantum theory of harmonic oscillators. On the negative real axis, $G(-\e^u)$ possesses a sign-constant asymptotic expansion in powers of $1/u^2$ for large positive $u$. We have shown that optimal truncation leaves exponentially suppressed terms of the form 
\[
\e^{-u}\sum_{k=0}^\infty \frac{P_k(x)}{u^k},
\]
with $P_0(x)=x-\tfrac23$ and $P_k(x)$ a polynomial of degree $2k+1$, evaluated at $x=u/2-\lfloor u/2\rfloor$. At large $k$, these polynomials closely approximate sinusoids. Their amplitude grows factorially with~$k$ and their phase increases linearly, with
\[
P_k(x)\sim \frac{1}{\sqrt{2\pi}} \sin((2k+1)C-2\pi x)R^{2k+1}\Gamma\!\left(k+\tfrac12\right),
\]
where $C$ and $R$ are asymptotic constants satisfying $R\exp(\im\,C)=\sqrt{-1/(2+\pi\im)}$.

We observe that $G(-\e^{u})$ is precisely a Fermi--Dirac integral of order $-\frac{3}{2}$, $G(-\e^{u})=-F_{-3/2}(u)$. The present analysis therefore extends directly to the full class of Fermi--Dirac integrals through~\cite[Eq.~(3.3)]{TMG}.

In this work, we were content to show that the sub-leading term $F_k(x,1)$, in \eqref{Pleading} for $P_k(x)$, is exponentially suppressed at large $k$, as shown in \eqref{Pkasy}. It is possible to go further and develop a complete asymptotic expansion for $F_k(x,1)$, in terms of gamma functions with descending arguments. Such behaviour is characteristic of late coefficients in resurgence theory (see, e.g.,~\cite{MBCH,MBDB}). The core approach would involve substituting truncated asymptotic expansions for the dominant terminant functions (expressed in ascending powers of $v$) into the integral representation \eqref{F1expr}, and subsequently estimating their error terms. We leave these details to the interested reader.

\appendix

\section{The terminant function}\label{terminant}

In this appendix, we summarise key properties of the terminant function $T_p(z)$, first introduced in~\cite{FWJO}. It is defined in terms of the incomplete gamma function by
\begin{equation}\label{Tpdef}
T_p (z) = \frac{\e^{\pi \im p} \Gamma (p)}{2\pi \im}\Gamma (1 - p,z),
\end{equation}
for $p \in \mathbb{C}$ and $z$ lying on the Riemann surface of the logarithm.

By taking complex conjugates of both sides of \eqref{Tpdef}, we arrive at the conjugacy relation
\begin{equation}\label{Tconj}
\overline {T_p (z)}  =  - \e^{ - 2\pi \im\bar p} T_{\bar p} (\bar z).
\end{equation}

The function $T_p(z)$ admits the integral representation~\cite[Eq.~(3.9)]{FWJO}
\begin{equation}\label{Tintegrals}
T_p(z) 
= \e^{\pi \im p}\frac{\e^{-z}}{2\pi \im} \int_0^{+\infty} \frac{ \e^{-zs}s^{p-1}}{1+s}\id s 
= \e^{\pi \im p}\frac{\e^{-z}}{2\pi \im} \int_{-\infty}^{+\infty} \frac{\e^{-z\e^t}\e^{pt} }{1+\e^t}\id t,
\end{equation}
valid for $\Re(p)>0$ and $\left|\arg z\right|<\tfrac{\pi}{2}$. By changing the integration variable in the first integral and applying analytic continuation in $z$, one obtains
\begin{equation}\label{Tintegral}
T_p (z) = \e^{\pi \im p} \frac{\e^{ - z} z^{ - p} }{2\pi \im}\int_0^{ + \infty } \frac{\e^{ - s} s^{p - 1} }{1 + s/z}\id s,
\end{equation}
valid for $\Re(p)>0$ and $\left|\arg z\right|<\pi$.

Setting $p = z + \alpha$ with bounded $\alpha \in \mathbb{C}$ in \eqref{Tintegrals}, an application of the saddle point method~\cite[\href{http://dlmf.nist.gov/2.4.iv}{\S2.4(iv)}]{DLMF} yields
\begin{equation}\label{asymp1}
T_p(z) 
= \e^{\pi \im p}\frac{\e^{-2z}}{2\pi \im} \int_{-\infty}^{+\infty} \e^{-z(\e^t-t-1)} \frac{\e^{\alpha t}}{1+\e^t}\id t  
\sim \e^{\pi \im p}\frac{\e^{-2z}}{2\pi \im} \sum_{k=0}^\infty \frac{A_k(\alpha)}{z^{k+1/2}},
\end{equation}
as $z \to \infty$ in the sector $\left|\arg z\right|\le\frac{\pi}{2}-\delta$ ($<\frac{\pi}{2}$), where
\begin{equation}\label{Acoeff}
A_k(\alpha) 
= \frac{\Gamma\!\left(k+\tfrac{1}{2}\right)}{2\pi \im} 
  \oint_{|t|=1} \frac{\e^{\alpha t}}{1+\e^t} 
  \frac{\d t}{\left(\e^t - t - 1\right)^{k+1/2}}.
\end{equation}
If, in addition, $\Re(\alpha) > 0$, then \eqref{asymp1} remains valid in the closed sector $\left|\arg z\right| \le \frac{\pi}{2}$. This can be seen by deforming the contour of integration to a path along which $\Re(z(\e^t-t-1))>0$. 

An asymptotic expansion valid for $\arg z = \pi$ can be obtained by suitably adapting the notation of~\cite[Eqs.~(5.2) and (4.1)]{FWJO2}. This yields
\begin{equation}\label{asymp2}
T_p\!\left(z\e^{\pi \im}\right) \sim \frac{1}{2} - \frac{\im}{2\pi} \sum_{k=0}^\infty \frac{B_k(\alpha)}{z^{k+1/2}},\quad p=z+\alpha,
\end{equation}
as $z\to+\infty$, where
\begin{align*}
B_k(\alpha) 
& = \frac{\Gamma\!\left(k+\tfrac{1}{2}\right)}{2\pi \im} 
  \oint_{|s-1|=\frac{1}{2}} \frac{s^{\alpha-1}}{1-s} 
  \frac{\d s}{\left(s-\log (s)-1\right)^{k+1/2}}
\\ & = \frac{\Gamma\!\left(k+\tfrac{1}{2}\right)}{2\pi \im} 
  \oint_{|t|=1} \frac{\e^{\alpha t}}{1-\e^t} 
  \frac{\d t}{\left(\e^t - t - 1\right)^{k+1/2}}.
\end{align*}
We remark that the coefficients $B_k(\alpha)$ are real when $\alpha$ is real. 

We will also use the following integral representation~\cite[p.~522, Eq.~(37.3.47)]{NT}:
\begin{equation}\label{Tpoleint}
T_p\! \left( z \e^{\pi \im} \right)
= \frac{\e^z}{2\pi \im} \cupint_0^{+\infty} \e^{-z s} \frac{s^{p-1}}{1 - s} \id s
= \frac{\e^z}{2\pi \im} \cupint_{-\infty}^{+\infty} \frac{\e^{-z \e^t} \e^{p t}}{1 - \e^t} \id t,
\end{equation}
which is valid for $\Re(p) > 0$ and $\left|\arg z\right| < \frac{\pi}{2}$. The integration contours are chosen to pass below the singularities at $s=1$ and $t=0$, respectively.

\section{Asymptotic behaviour of the coefficients \texorpdfstring{$\Delta_k(x)$}{Delta\_k(x)}}\label{Dkasyproof}

In this appendix, we derive the asymptotic formula \eqref{Dkasy} for the polynomials $\Delta_k(x)$.

We begin by applying the change of integration variables $s = u\e^t$ to the integral definition of the gamma function. Provided that $\Re(u-2x+\tfrac{1}{2})>0$, this yields
\[
\Gamma\!\left(u-2x+\tfrac{1}{2}\right)
 = \int_0^{+\infty} \e^{-s}s^{u-2x-1/2}\id s
 = u^{u-2x+1/2}\e^{-u}
   \int_{-\infty}^{+\infty}
   \exp\!\big(-u(\e^t-t-1)\big)\,
   \e^{(1/2-2x)t}\id t.
\]
Applying the saddle point method~\cite[\href{http://dlmf.nist.gov/2.4.iv}{\S2.4(iv)}]{DLMF}, we obtain the asymptotic expansion as $u\to+\infty$:
\[
\frac{\e^u}{\sqrt{2\pi}}\,\frac{\Gamma(u-2x+\tfrac12)}{u^{u-2x}} = 
\sqrt {\frac{u}{2\pi }}\int_{-\infty}^{+\infty}
   \exp\!\big(-u(\e^t-t-1)\big)\,
   \e^{(1/2-2x)t}\id t \sim \sum_{k=0}^\infty\frac{\Delta_k(x)}{u^k},
\]
where, for any real $x$, the coefficients are given by
\[
\Delta_k(x)=\frac{\Gamma\!\left(k+\tfrac12\right)}
     {\sqrt{2\pi}}
\frac{1}{2\pi \im}
\oint_{|t|=1}
\frac{  \e^{(1/2-2x)t}}{\left(\e^t-t-1\right)^{k+1/2}
}\id t.
\]
Under the assumption that $-\frac{k}{2}<x<\frac{1}{4}$ and $k\ge 1$, we can deform the circular integration contour onto the boundary of the domain $D$ shown in Figure~\ref{fig1}. We thus obtain
\begin{align*}
\Delta_k(x)
&= \frac{\Gamma\!\left(k+\tfrac12\right)}
     {\sqrt{2\pi}}
\frac{1}{2\pi \im}
   \sum_{m=\pm1}
   \int_{\mathscr{C}^{(m)}}
   \frac{\e^{(1/2-2x)t}}{\left(\e^t-t-1\right)^{k+1/2}}\id t \\
&= \frac{\Gamma\!\left(k+\tfrac12\right)}
     {\sqrt{2\pi}}
\frac{1}{2\pi \im}
   \sum_{m=\pm1}
   \frac{1}{\left(-2\pi \im m\right)^{k+1/2}}
   \int_{\mathscr{C}^{(m)}}
   \exp\!\left(-\left(k+\tfrac{1}{2}\right)
   \log\!\left(\frac{\e^t-t-1}{-2\pi \im m}\right)\right)
   \e^{(1/2-2x)t}\id t.
\end{align*}
By analytic continuation, this representation remains valid in the larger domain $x>-\frac{k}{2}$. Next, we apply the saddle point method to each of these two integrals. By carefully choosing the branches of the square roots arising during the procedure (following the prescription given for the specific case $x=\frac{1}{4}$ in~\cite[\S3.3]{WB2}), we deduce that
\begin{align*}
\Delta_k(x)&=\frac{\Gamma\!\left(k+\frac12\right)}
     {\sqrt{k+\frac12}}\frac{1}{2\pi \im}
\sum_{m=\pm1}m \frac{\e^{2\pi \im(1/2-2x)m}}
        {\left(-2\pi \im m\right)^k}
\left(1+\mathcal{O}\!\left(\frac1k\right)\right) \\
&=\frac{2\,\Gamma\!\left(k+\frac12\right)}
     {\sqrt{k+\frac12}\left(2\pi\right)^{k+1}}
\left(\sin\!\left(4\pi x-\tfrac{\pi}{2}k\right)
+\mathcal{O}\!\left(\frac1k\right)
\right),
\end{align*}
as $k\to+\infty$ for any fixed real $x$. Finally, substituting the asymptotic behaviour 
\[
\frac{\Gamma\!\left( k + \frac{1}{2} \right)}{\sqrt {k + \frac{1}{2}} } = \Gamma (k)\!\left( 1 + \mathcal{O}\!\left( {\frac{1}{k}} \right) \right),\quad k\to+\infty,
\]
(cf.~\cite[Eq.~\href{http://dlmf.nist.gov/5.11.E13}{(5.11.13)}]{DLMF}) yields the desired formula \eqref{Dkasy}.

\section*{Acknowledgments}

DB is grateful for discussions with David Gross and Gerard {}'t Hooft, at {\em Bohr-100: Current Themes in Theoretical Physics}, held at the Niels Bohr Institute, in 2022, and with Ovidiu Costin, Daniele Dorigoni and Gerald Dunne, at {\em Resurgence and Modularity in QFT and String Theory}, held at the Galileo Galilei Institute, in 2024. The authors gratefully acknowledge the referees for their insightful comments and valuable suggestions, which have significantly improved the clarity and presentation of the paper.

\section*{Declarations}

\paragraph{\textbf{Conflict of interest}} On behalf of all authors, the corresponding author states that there is no conflict of interest.

\section*{Data availability}

The authors confirm that the data supporting the findings of this study are available within the article. No additional datasets were generated or analysed during the current study.

\let\oldbibitem\bibitem
\def\bibitem{\par\addvspace{3pt}\oldbibitem}


\begin{thebibliography}{99}

\bibitem{MVB}
M.~V.~Berry, Infinitely many Stokes smoothings in the gamma function, \emph{Proc. Roy. Soc. London Ser. A} \textbf{434} (1991), pp.~465--472.

\bibitem{MBCH}
M.~V.~Berry and C.~J.~Howls, Hyperasymptotics for integrals with saddles, \emph{Proc. Roy. Soc. London Ser. A} \textbf{434} (1991), pp.~657--675.

\bibitem{MBDB}
M. Borinsky and D. Broadhurst, Resonant resurgent asymptotics from quantum field
theory, \emph{Nucl. Phys. B} \textbf{981} (2022), Article 115861, \url{https://arxiv.org/abs/2202.01513}

\bibitem{WB1}
W.~G.~C.~Boyd, Error bounds for the method of steepest descents, \emph{Proc. Roy. Soc. London Ser. A} \textbf{440} (1993),
pp.~493--518.

\bibitem{WB}
W.~G.~C.~Boyd, Gamma function asymptotics by an extension of the method of steepest descents, \emph{Proc. Roy. Soc. London Ser. A} \textbf{447} (1994), pp.~609--630.

\bibitem{WB2}
W.~G.~C.~Boyd, Approximations for the late coefficients in asymptotic expansions arising in the method of steepest descents, \emph{Methods Appl. Anal.} \textbf{2} (1995), pp.~475--489.

\bibitem{OC}
O.~Costin and S.~Garoufalidis, Resurgence of the fractional polylogarithms, \emph{Math. Res. Lett.} \textbf{16} (2009), pp.~817--826, \url{https://arxiv.org/abs/math/0701743}

\bibitem{TMG}
T.~M.~Garoni, N.~E.~Frankel, and M.~L.~Glasser, Complete asymptotic expansions of the Fermi--Dirac integrals, \emph{J. Math. Phys.} \textbf{42} (2001), pp.~1860--1868.

\bibitem{GH}
G.~'t~Hooft, The hidden ontological variable in quantum harmonic oscillators, \emph{Front. Quantum Sci. Technol.} \textbf{3} (2024), Article 1505593, \url{https://arxiv.org/abs/2407.18153}

\bibitem{FJ}
F.~Johansson, Rigorous high-precision computation of the Hurwitz zeta function and its derivatives, \emph{Numer. Algor.} \textbf{69} (2015), pp.~253--270, \url{https://arxiv.org/abs/1309.2877}

\bibitem{GN}
G.~Nemes, Error bounds and exponential improvements for the asymptotic expansions of the gamma function and its reciprocal, \emph{Proc. Roy. Soc. Edinburgh Sect. A} \textbf{145} (2015), pp.~571--596, \url{https://arxiv.org/abs/1310.0166}

\bibitem{DLMF}
\emph{NIST Digital Library of Mathematical Functions}. \url{https://dlmf.nist.gov/}, Release 1.2.6 of 2026-03-15. F.~W.~J.~Olver, A.~B.~Olde Daalhuis, D.~W.~Lozier, B.~I.~Schneider, R.~F.~Boisvert, C.~W.~Clark, B.~R.~Miller, B.~V.~Saunders, H.~S.~Cohl, and M.~A.~McClain, eds.

\bibitem{ODO}
A.~B.~Olde Daalhuis and F.~W.~J.~Olver, Exponentially improved asymptotic solutions of
ordinary differential equations. II. Irregular singularities of rank one, \emph{Proc. Roy. Soc. London Ser. A} \textbf{445} (1994),
pp.~39--56.

\bibitem{FWJO}
F.~W.~J.~Olver, On Stokes' phenomenon and converging factors, in: R.~Wong, Ed., \emph{Asymptotic and Computational Analysis}, Marcel Dekker, New York, 1990, pp. 329--355. 

\bibitem{FWJO2}
F.~W.~J.~Olver, Uniform, exponentially improved, asymptotic expansions for the generalized exponential integral, \emph{SIAM J. Math. Anal.} \textbf{22} (1991), pp.~1460--1474.

\bibitem{Pari}
The PARI Group, University of Bordeaux, \emph{Pari/GP}, version 2.17.2 (2025), \url{http://pari.math.u-bordeaux.fr}

\bibitem{NT}
N.~M.~Temme, \emph{Asymptotic Methods for Integrals}, World Scientific, London, 2015.

\bibitem{DW}
D.~C.~Wood, The computation of polylogarithms, Technical Report 15-92, University of Kent (1992), \url{https://www.cs.kent.ac.uk/pubs/1992/110/content.pdf}

\end{thebibliography}
\end{document}